\documentclass{gtart}
\input gtmonout
\volumenumber{1}
\volumeyear{1998}
\volumename{The Epstein birthday schrift}
\pagenumbers{23}{49}
\received{15 November 1997}
\published{21 October 1998}
\papernumber{2}

\newtheorem{conjecture}{Conjecture}[section]

\newtheorem{lemma}[conjecture]{Lemma}
\newtheorem{proposition}[conjecture]{Proposition}

\newtheorem{theorem}[conjecture]{Theorem}

\begin{document}

\title{A brief survey of the deformation theory\\of Kleinian groups} 
\shorttitle{Deformation theory of Kleinian groups}
\author{James W Anderson} 

\address{Faculty of Mathematical
Studies\\University of Southampton\\
Southampton, SO17 1BJ, England}
\email{jwa@maths.soton.ac.uk}

\begin{abstract} 
We give a brief overview of the current state of the study of the
deformation theory of Kleinian groups.  The topics covered include the
definition of the deformation space of a Kleinian group and of several
important subspaces; a discussion of the parametrization by
topological data of the components of the closure of the deformation
space; the relationship between algebraic and geometric limits of
sequences of Kleinian groups; and the behavior of several
geometrically and analytically interesting functions on the
deformation space.
\end{abstract}



\makeshorttitle

\cl{{\small\em Dedicated to David Epstein on the occasion of his 60th
birthday}}

\section{Introduction}
\label{introduction}

Kleinian groups, which are the discrete groups of orientation
preserving isometries of hyperbolic space, have been studied for a
number of years, and have been of particular interest since the work
of Thurston in the late 1970s on the geometrization of compact
$3$--manifolds.  A Kleinian group can be viewed either as an isolated,
single group, or as one of a member of a family or continuum of
groups.  

In this note, we concentrate our attention on the latter scenario,
which is the deformation theory of the title, and attempt to give a
description of various of the more common families of Kleinian groups
which are considered when doing deformation theory.  No proofs are
given, though it is hoped that reasonable coverage of the current
state of the subject is given, and that ample references have been
given for the interested reader to venture boldly forth into the
literature.

It is possible to consider the questions raised here in much more
general settings, for example for Kleinian groups in $n$--dimensions
for general $n$, but that is beyond the scope of what is attempted
here.  Some material on this aspect of the question can be found in
Bowditch \cite{bowditch-higher} and the references contained therein.

The author would like to thank Dick Canary, Ed Taylor, and Brian
Bowditch for useful conversations during the preparation of this
work, as well as the referee for useful comments.

\section{The deformation spaces}
\label{deformation-spaces}

We begin by giving a few basic definitions of the objects considered
in this note, namely Kleinian groups.  We go on to define and describe
the basic structure of the deformation spaces we are considering
herein.

A {\em Kleinian group} is a discrete subgroup of ${\rm PSL}_2({\bf C})
={\rm SL}_2({\bf C})/\{ \pm {\rm I}\}$, which we view as acting both
on the Riemann sphere $\overline{{\bf C}}$ by M\"obius transformations
and on real hyperbolic $3$--space ${\bf H}^3$ by isometries, where the
two actions are linked by the Poincar\'e extension. 

The action of an infinite Kleinian group $\Gamma$ partitions
$\overline{{\bf C}}$ into two sets, the {\em domain of discontinuity}
$\Omega(\Gamma)$, which is the largest open subset of $\overline{{\bf
C}}$  on which $\Gamma$ acts discontinuously, and the {\em limit set}
$\Lambda(\Gamma)$.  If $\Lambda(\Gamma)$ contains two or fewer points,
$\Gamma$ is  {\em elementary}, otherwise $\Gamma$ is {\em
non-elementary}.  For a non-elementary Kleinian group $\Gamma$, the
limit set $\Lambda(\Gamma)$ can also be described as the smallest
non-empty closed subset of $\overline{{\bf C}}$ invariant under
$\Gamma$.  We refer the reader to Maskit \cite{maskit-book} or
Matsuzaki and Taniguchi \cite{matsuzaki-taniguchi} as a
reference for the basics of Kleinian groups.

An isomorphism $\varphi\co  \Gamma\rightarrow\Phi$ between Kleinian
groups $\Gamma$ and $\Phi$ is {\em type-preserving} if, for
$\gamma\in\Gamma$, we have that $\gamma$ is parabolic if and only if
$\varphi(\gamma)$ is parabolic.

A Kleinian group is {\em convex cocompact} if its convex core is
compact; recall that the {\em convex core} associated to a Kleinian
group $\Gamma$ is the minimal convex submanifold of ${\bf H}^3/\Gamma$
whose inclusion is a homotopy equivalence.  More generally, a Kleinian
group is {\em geometrically finite} if it is finitely generated and if
its convex core has finite volume.  This is one of several equivalent
definitions of geometric finiteness; the interested reader is referred
to Bowditch \cite{bowditch} for a complete discussion.

A Kleinian group $\Gamma$ is {\em topologically tame} if its
corresponding quotient $3$--man\-ifold ${\bf H}^3/\Gamma$ is homeomorphic
to the interior of a compact $3$--manifold.  Geometrically finite
Kleinian groups are topologically tame.  It was conjectured by Marden
\cite{marden} that all finitely generated Kleinian groups are
topologically tame.

A compact $3$--manifold $M$ is {\em hyperbolizable} if there exists a
hyperbolic $3$--manifold $N ={\bf H}^3/\Gamma$ homeomorphic to the
interior of $M$.  Note that a hyperbolizable $3$--manifold $M$ is 
necessarily orientable; {\em irreducible}, in that every embedded
$2$--sphere in $M$ bounds a $3$--ball in $M$; and {\em atoroidal}, in
that every embedded torus $T$ in $M$ is homotopic into $\partial M$.
Also, since the universal cover ${\bf H}^3$ of $N$ is contractible,
the fundamental group of $M$ is isomorphic to $\Gamma$.  For a
discussion of the basic theory of $3$--manifolds, we refer the reader
to Hempel \cite{hempel} and Jaco \cite{jaco}.

Keeping to our viewpoint of a Kleinian group as a member of a family
of groups, throughout this survey we view a Kleinian group as the
image $\rho(G)$ of a representation $\rho$ of a group $G$ into ${\rm 
PSL}_2({\bf C})$.  Unless explicitly stated otherwise, we assume that
$G$ is finitely generated, torsion-free, and non-abelian, so that
in particular $\rho(G)$ is non-elementary.

\subsection{The representation varieties ${\cal HOM}(G)$ and\nl ${\rm R}(G)
={\cal HOM}(G)/{\rm PSL}_2({\bf C})$}

The most basic of the deformation spaces is the {\em representation
variety} ${\cal HOM}(G)$ which is the space of all representations of
$G$ into ${\rm PSL}_2({\bf C})$ with the following topology.  Given a
set of generators $\{ g_1,\ldots, g_k\}$ for $G$, we may naturally
view ${\cal HOM}(G)$ as a subset of ${\rm PSL}_2({\bf C})^k$, where a
representation $\rho\in{\cal HOM}(G)$ corresponds to the $k$--tuple 
$(\rho(g_1),\ldots,\rho(g_k))$ in ${\rm PSL}_2({\bf C})^k$.  The
defining polynomials of this variety are determined by the relations
in $G$.  In particular, if $G$ is free, then ${\cal HOM}(G) ={\rm
PSL}_2({\bf C})^k$.  It is easy to see that ${\cal HOM}(G)$ is a
closed subset of ${\rm PSL}_2({\bf C})^k$.  

The representations in ${\cal HOM}(G)$ are {\em unnormalized}, in the
sense that there is a natural free action of ${\rm PSL}_2({\bf C})$ on
${\cal HOM}(G)$ by conjugation.  Depending on the particular question
being addressed, it is sometimes preferable to remove the ambiguity of
this action and form the quotient space ${\rm R}(G)={\cal HOM}(G)/{\rm
PSL}_2({\bf C})$.

Though a detailed description is beyond the scope of this survey, we
pause here to mention work of Culler and Shalen \cite{cs-varieties},
\cite{cs-bounded}, in which a slight variant of the representation
variety as defined above plays a fundamental role, and which has
inspired further work of Morgan and Shalen \cite{ms-I}, 
\cite{ms-II}, \cite{ms-III} and Culler, Gordon, Luecke, and Shalen
\cite{cgls}.  The basic object here is not the space ${\rm R}(G)$ of 
all representations of $G$ into ${\rm PSL}_2({\bf C})$ as defined
above, but instead the related space ${\rm X}(G)$ of all
representations of $G$ into ${\rm SL}_2({\bf C})$, modulo the action
of ${\rm SL}_2({\bf C})$.  The introduction of this space ${\rm X}(G)$
does beg the question of when a representation of $G$ into ${\rm
PSL}_2({\bf C})$ can be lifted to a representation of $G$ into ${\rm
SL}_2({\bf C})$.  We note in passing that this question of lifting
representations has been considered by a number of authors, including
Culler, Kra, and Thurston, to name but a few; we refer the reader to
the article by Kra \cite{kra-lifting} for exact statements and a
review of the history, including references.

By considering the global structure of the variety ${\rm X}(G)$ in the
case that $G$ is the fundamental group of a compact, hyperbolizable
$3$--manifold $M$, and in particular the ideal points of its
compactification, Culler and Shalen \cite{cs-varieties} are able to
analyze the actions of $G$ on trees, which in turn has connections
with the existence of essential incompressible surfaces in $M$, finite
group actions on $M$, and has particular consequences in the case
that $M$ is the complement of a knot in ${\bf S}^3$. We refer the
reader to the excellent survey article by Shalen \cite{shalen-survey},
as well as to the papers cited above.

\subsection{The spaces ${\cal HOM}_{{\rm T}}(G)$ and ${\rm R}_{{\rm
T}}(G) ={\cal HOM}_{{\rm T}}(G)/{\rm PSL}_2({\bf C})$ of the minimally
parabolic representations}

Let ${\cal HOM}_{{\rm T}}(G)$ denote the subspace of ${\cal HOM}(G)$
consisting of those representations $\rho$ for which $\rho(g)$ is
parabolic if and only if $g$ lies in a rank two free abelian subgroup
of $G$.  We refer to ${\cal HOM}_{{\rm T}}(G)$ as the space of {\em
minimally parabolic} representations of $G$. In particular, if $G$
contains no ${\bf Z}\oplus {\bf Z}$ subgroups, then the image
$\rho(G)$ of every $\rho$ in ${\cal HOM}_{{\rm T}}(G)$ is {\em purely
loxodromic}, in that every non-trivial element of $\rho(G)$ is
loxodromic.  Set ${\rm R}_{{\rm T}}(G)={\cal HOM}_{{\rm T}}(G)/{\rm
PSL}_2({\bf C})$. 

\subsection{The spaces ${\cal D}(G)$ and ${\rm AH}(G) ={\cal
D}(G)/{\rm PSL}_2({\bf C})$ of discrete, faithful representations} 

Let ${\cal D}(G)$ denote the subspace of ${\cal HOM}(G)$ consisting of
the discrete, faithful representations of $G$, that is, the injective
homomorphisms of $G$ into ${\rm PSL}_2({\bf C})$ with discrete image.
For the purposes of this note, the space ${\cal D}(G)$ is our
universe, as it is the space of all Kleinian groups isomorphic to
$G$.  Set ${\rm AH}(G)={\cal D}(G)/{\rm PSL}_2({\bf C})$.  

We note that there exists an equivalent formulation of ${\rm AH}(G)$
in terms of manifolds.  Given a hyperbolic $3$--manifold $N$, let ${\rm
H}(N)$ denote the set of all pairs $(f,K)$, where $K$ is a hyperbolic
$3$--manifold and $f\co N\rightarrow K$ is a homotopy equivalence, modulo
the equivalence relation $(f,K)\sim (g,L)$ if there exists an
orientation preserving isometry $\alpha\co  K\rightarrow L$ so that
$\alpha\circ f$ is homotopic to $g$.  The topology on ${\rm H}(N)$ is
given by noting that, if we let $\Gamma\subset {\rm PSL}_2({\bf C})$
be a choice of conjugacy class of the fundamental group of $N$, then
each element $(f,K)$ in ${\rm H}(N)$ gives rise to a discrete,
faithful representation $\varphi =f_\ast$ of $\Gamma$ into ${\rm
PSL}_2({\bf C})$, with equivalent points in ${\rm H}(N)$ giving rise
to conjugate representations into ${\rm PSL}_2({\bf C})$.
Hence, equipping ${\rm H}(N)$ with this topology once again gives rise
to ${\rm AH}(G)$ with $G =\pi_1(N)$.

The following theorem, due to J\o rgensen, describes the fundamental
property of ${\cal D}(G)$, namely that the limit of a sequence of
elements of ${\cal D}(G)$ is again an element of ${\cal D}(G)$. 

\begin{theorem}[J\o rgensen \cite{jorgensen}] ${\cal D}(G)$ is a
closed subset of ${\cal HOM}(G)$.
\label{deformation-space-closed}
\end{theorem}

There is one notable case in which ${\rm AH}(G)$ is completely
understood, namely in the case that $G$ is the fundamental group of a
compact, hyperbolizable $3$--manifold $M$ whose boundary is the union
of a (possibly empty) collection of tori.  In this case, the
hyperbolic structure on the interior of $M$ is unique, by the
classical Rigidity Theorem of Mostow, for closed manifolds, and
Prasad, for manifolds with non-empty toroidal boundary.  Rephrasing
this statement as a statement about deformation spaces yields the
following. 

\begin{theorem}[Mostow \cite{mostow} and Prasad \cite{prasad}]
Suppose that $G$ is the fundamental group of a compact, orientable
$3$--manifold $M$ whose boundary is the union of a (possibly empty)
collection of tori.  Then, ${\rm AH}(G)$ either is empty or consists
of a single point.
\label{mostow-rigidity}
\end{theorem}

Given this result, it will cause us no loss of generality to assume
that henceforth all Kleinian groups have infinite volume quotients.

\subsection{The spaces ${\cal P}(G)$ and ${\rm MP}(G) ={\cal
P}(G)/{\rm PSL}_2({\bf C})$ of geometrically finite, minimally
parabolic representations}

Let ${\cal P}(G)$ denote the subset of ${\cal D}(G)$ consisting of
those representations $\rho$ with geometrically finite, minimally
parabolic image $\rho(G)$.  In particular, if $G$ contains no ${\bf
Z}\oplus {\bf Z}$ subgroups, then the image $\rho(G)$ of every
representation $\rho\in {\cal P}(G)$ is convex cocompact.  Set ${\rm
MP}(G) ={\cal P}(G)/{\rm PSL}_2({\bf C})$, and note that since ${\rm
PSL}_2({\bf C})$ is connected, the quotient map gives a one-to-one
correspondence between the connected components of ${\cal P}(G)$ and
those of ${\rm MP}(G)$. 

It is an immediate consequence of the Core Theorem of Scott
\cite{scott} and the Hyperbolization Theorem of Thurston that if ${\cal
D}(G)$ is non-empty, then ${\cal P}(G)$ is non-empty.  For a
discussion of the Hyperbolization Theorem, see Morgan \cite{morgan},
Otal and Paulin \cite{otal-paulin}, or Otal \cite{otal-fibered} for
the fibered case.

We note here that, if there exists a geometrically finite, minimally
parabolic representation of $G$ into ${\rm PSL}_2({\bf C})$, then
in general there exist many geometrically finite representations which
are not minimally parabolic, which can be constructed as limits of
the geometrically finite, minimally parabolic representations.  This
construction has been explored in detail for a number of cases by
Maskit \cite{maskit-squeezing} and Ohshika \cite{ohshika-squeezing}.

In the case that $G$ is itself a geometrically finite, minimally
parabolic Kleinian group, the structure of ${\rm MP}(G)$ is
fairly well understood, both as a subset of ${\rm AH}(G)$ and in terms
of how the components of ${\rm MP}(G)$ are parametrized by
topological data.  We spend the remainder of this section making these
statements precise.

We begin with the Quasiconformal Stability Theorem of Marden
\cite{marden}.

\begin{theorem}[Marden \cite{marden}] If $G$ is a geometrically
finite, minimally para\-bolic Kleinian group, then ${\rm MP}(G)$ is an
open subset of ${\rm R}(G)$. 
\label{quasiconformal-stability}
\end{theorem}

As a converse to this, we have the Structural Stability Theorem of
Sullivan \cite{sullivan-stability}.  We note here that the versions of
the Theorems of Marden and Sullivan given here are not the strongest,
but are adapted to the point of view taken in this paper. The general
statements holds valid in slices of ${\rm AH}(G)$ in which a certain 
collection of elements of $G$ are required to have parabolic image,
not just those which belong to ${\bf Z}\oplus {\bf Z}$ subgroups.

\begin{theorem}[Sullivan \cite{sullivan-stability}] Let $G$ be a
finitely generated, torsion-free, non-elementary Kleinian group.  If
there exists an open neighborhood of the identity representation in
${\rm R}(G)$ which lies in ${\rm AH}(G)$, then $G$ is geometrically
finite and minimally parabolic.  
\label{structural-stability}
\end{theorem}

Combining these, we see that ${\rm MP}(G)$ is the interior of ${\rm
AH}(G)$.  A natural question which arises from this is whether there
are points of ${\rm AH}(G)$ which do not lie in the closure of ${\rm
MP}(G)$.

\begin{conjecture}[Density conjecture] ${\rm AH}(G)$ is the closure
of ${\rm MP}(G)$.
\label{density-conjecture}
\end{conjecture}

This Conjecture is due originally to Bers in the case that $G$ is the
fundamental group of a surface, see Bers \cite{bers-boundaries}, and 
extended by Thurston to general $G$.

There has been a good deal of work in the past couple of years on the
global structure of ${\rm MP}(G)$ and its closure.  We begin with an
example to show that there exist groups $G$ for which ${\rm MP}(G)$
is disconnected; the example we give here, in which ${\rm MP}(G)$ has 
finitely many components, comes from the discussion in Anderson and
Canary \cite{ac-books}.

Let $T$ be a solid torus and for large $k$, let $A_1,\ldots, A_k$ be
disjoint embedded annuli in $\partial T$ whose inclusion into $T$
induces an isomorphism of fundamental groups.  For each $1\le j\le k$,
let $S_j$ be a compact, orientable surface of genus $j$ with a single
boundary component, and let $Y_j =S_j\times I$, where $I$ is a closed
interval.  Construct a compact $3$--manifold $M$ by attaching the
annulus $\partial S_j\times I$ in $\partial Y_j$ to the annulus $A_j$
in $\partial T$.  The resulting $3$--manifold $M$ is compact and
hyperbolizable $3$--manifold and has fundamental group $G$.  This
$3$--manifold is an example of a {\em book of I--bundles}.  Let $\rho$
be an element of ${\rm MP}(G)$ for which the interior of $M$ is
homeomorphic to ${\bf H}^3/\rho(G)$.

Let $\tau$ be a permutation of $\{ 1,\ldots, k\}$, and consider now
the manifold $M_\tau$ obtained by attaching the annulus $\partial
S_j\times I$ in $\partial Y_j$ to the annulus $A_{\tau(j)}$ in
$\partial T$.  By construction, $M_\tau$ is compact and
hyperbolizable, and has fundamental group $G$; let $\rho_\tau$ be an
element of ${\rm MP}(G)$ for which the interior of $M_\tau$ is
homeomorphic to ${\bf H}^3/\rho_\tau(G)$.  Since $M$ and $M_\tau$ have
isomorphic fundamental groups, they are homotopy equivalent. However,
in the case that $\tau$ is not some power of the cycle $(12\cdots k)$,
then there does not exist an orientation preserving homeomorphism
between $M$ and $M_\tau$, and hence $\rho$ and $\rho_\tau$ lie in
different components of ${\rm MP}(G)$.

In the general case that $G$ is finitely generated and does not split
as a free product, there exists a characterization of the components
of both ${\rm MP}(G)$ and its closure $\overline{{\rm MP}(G)}$ in
terms of the topology of a compact, hyperbolizable $3$--manifold $M$
with fundamental group $G$.  This characterization combines work of
Canary and McCullough \cite{canary-mccullough} and of Anderson,
Canary, and McCullough \cite{acm}.  We need to develop a bit of
topological machinery before discussing this characterization.

For a compact, oriented, hyperbolizable $3$--manifold $M$ with
non-empty, incompressible boundary, let ${\cal A}(M)$ denote the set
of {\em marked homeomorphism types of compact, oriented $3$--manifolds
homotopy equivalent to $M$}. Explicitly, ${\cal A}(M)$ is the set of
equivalence classes of pairs $(M',h')$, where $M'$ is a compact,
oriented, irreducible $3$--manifold and $h'\co M\to M'$ is a homotopy
equivalence, and where two pairs $(M_1,h_1)$ and $(M_2,h_2)$ are
equivalent if there exists an orientation preserving homeomorphism
$j\co M_1\to M_2$ such that $j\circ h_1$ is homotopic to 
$h_2$.  Denote the class of $(M',h')$ in ${\cal A}(M)$ by
$[(M',h')]$. 

There exists a natural map $\Theta\co  {\rm AH}(\pi_1(M))\to {\cal
A}(M)$, defined as follows.  For $\rho\in {\rm AH}(\pi_1(M))$, let
$M_\rho$ be a compact core for $N_\rho ={\bf H}^3/\rho(\pi_1(M))$ and
let $r_\rho\co  M\to M_\rho$ be a homotopy equivalence such that 
$(r_\rho)_*\co \pi_1(M)\to\pi_1(M_\rho)$ is equal to $\rho$.  Set
$\Theta(\rho) =[(M_\rho, h_\rho)]$.  It is known that the restriction
of $\Theta$ to ${\rm MP}(\pi_1(M))$ is surjective, and that two
elements $\rho_1$ and $\rho_2$ of ${\rm  MP}(\pi_1(M))$ lie in the
same component of ${\rm MP}(\pi_1(M))$ if and only if
$\Theta(\rho_1)=\Theta(\rho_2)$.  Hence, $\Theta$ induces a one-to-one
correspondence between the components of ${\rm MP}(\pi_1(M))$ and the
elements of ${\cal A}(M)$; the reader is directed to Canary and
McCullough \cite{canary-mccullough} for complete details. 

Given a pair $M_1$ and $M_2$ of compact, hyperbolizable $3$--manifolds
with non-empty, incompressible boundary, say that a homotopy
equivalence $h\co M_1\to M_2$ is a {\em primitive shuffle} if there
exists a finite collection ${\cal V}_1$ of primitive solid torus
components of the characteristic submanifold $\Sigma(M_1)$ and a
finite collection ${\cal V}_2$ of solid torus components of
$\Sigma(M_2)$, so that $h^{-1}({\cal V}_2)={\cal V}_1$ and so that $h$
restricts to an orientation preserving homeomorphism from
$\overline{M_1-{\cal V}_1}$ to $\overline{M_2-{\cal V}_2}$; we do not
define the {\em characteristic submanifold} here, but instead refer
the reader to Canary and McCullough \cite{canary-mccullough}, Jaco and
Shalen \cite{jaco-shalen}, or Johannson \cite{johannson}.

Let $[(M_1,h_1)]$ and $[(M_2,h_2)]$ be two elements of ${\cal A}(M)$.
Say that $[(M_2,h_2)]$ is {\em primitive shuffle equivalent} to
$[(M_1,h_1)]$ if there exists a primitive shuffle $\varphi\co M_1\to M_2$
such that $[(M_2,h_2)] =[(M_2, \varphi\circ h_1)]$.  We note that when
$M$ is hyperbolizable, this gives an equivalence relation on ${\cal
A}(M)$, where each equivalence class contains finitely many elements
of ${\cal A}(M)$; let $\widehat{{\cal A}}(M)$ denote the set of
equivalence classes.  By considering the composition $\widehat\Theta
=q\circ\Theta$ of $\Theta$ with the quotient map $q\co  {\cal
A}(M)\rightarrow\widehat{{\cal A}}(M)$, we 
obtain the following complete enumeration of the components of
$\overline{{\rm MP}(\pi_1(M))}$. 

\begin{theorem}[Anderson, Canary, and McCullough \cite{acm}] Let $M$
be a compact, hyperbolizable $3$--manifold with non-empty,
incompressible boundary, and let $[(M_1,h_1)]$ and $[(M_2,h_2)]$ be
two elements of ${\cal A}(M)$.  The associated components of ${\rm
MP}(\pi_1(M))$ have intersecting closures if and only if $[(M_2,h_2)]$
is primitive shuffle equivalent to $[(M_1,h_1)]$.  In particular,
$\widehat\Theta$ gives a one-to-one correspondence between the
components of $\overline{{\rm MP}(\pi_1(M))}$ and the elements of
$\widehat{\cal A}$.
\label{enumeration}
\end{theorem}

Before closing this section, we highlight two consequences of the
analysis involved in the proof of Theorem \ref{enumeration}.  The
first involves the accumulation, or, more precisely, the lack thereof,
of components of ${\rm MP}(\pi_1(M))$.

\begin{proposition}{\rm(Anderson, Canary, and McCullough \cite{acm})}\stdspace
Let $M$ be a\break compact, hyperbolizable $3$--manifold with non-empty,
incompressible boundary.\break  Then, the components of ${\rm MP}(\pi_1(M))$
cannot accumulate in ${\rm AH}(\pi_1(M))$. In particular, the closure
$\overline{{\rm MP}(\pi_1(M))}$ of ${\rm MP}(\pi_1(M))$ is the union
of the closures of the components of ${\rm MP}(\pi_1(M))$.
\end{proposition}

The second involves giving a complete characterization, in terms of
the topology of $M$, as to precisely when $\overline{{\rm
MP}(\pi_1(M)))}$ has infinitely many components.  Recall that a
compact, hyperbolizable $3$--manifold $M$ with non-empty,
incompressible boundary has {\em double trouble} if there exists a
toroidal component $T$ of $\partial M$ and homotopically non-trivial
simple closed curves $C_1$ in $T$ and $C_2$ and $C_3$ in $\partial
M-T$ such that $C_2$ and $C_3$ are not homotopic in $\partial M$, but
$C_1$, $C_2$ and $C_3$ are homotopic in $M$.

\begin{theorem}[Anderson, Canary, and McCullough \cite{acm}] Let $M$
be a compact, hyperbolizable $3$--manifold with non-empty,
incompressible boundary.\break Then, $\overline{{\rm MP}(\pi_1(M))}$ has
infinitely many components if and only if $M$ has double trouble.
Moreover, if $M$ has double trouble, then ${\rm AH}(\pi_1(M))$ has
infinitely many components.
\label{infinitely-many-components}
\end{theorem}

\subsection{The spaces ${\cal QC}(G)$ and ${\rm QC}(G) ={\cal
QC}(G)/{\rm PSL}_2({\bf C})$ of quasiconformal deformations} 

In the case that $G$ is itself a finitely generated Kleinian group,
the classical deformation theory of $G$ consists largely of the study
of the space of {\em quasiconformal deformations} of $G$, which
consists of those representations of $G$ into ${\rm PSL}_2({\bf C})$
which are induced by a quasiconformal homeomorphism of the Riemann
sphere $\overline{\bf C}$. 

We do not give a precise definition here, but roughly, a {\em
quasiconformal homeomorphism} $\omega$ of $\overline{\bf C}$ is a
homeomorphism which distorts the standard complex structure on
$\overline{\bf C}$ by a bounded amount; the interested reader is
referred to Ahlfors \cite{ahlfors-quasiconformal} or to Lehto and
Virtanen \cite{lehto-virtanen} for a thorough
discussion of quasiconformality.  We do note that a quasiconformal
homeomorphism $\omega\co  \overline{{\bf C}}\rightarrow\overline{{\bf
C}}$ is completely determined (up to post-composition by a M\"obius
transformation) by the measurable function $\mu
=\omega_{\overline{z}}/\omega_z$, and that to every measurable
function $\mu$ on $\overline{\bf C}$ with
$\parallel\mu\parallel_\infty <1$ there exists a quasiconformal
homeomorphism $\omega$ of $\overline{\bf C}$ which solves the Beltrami
equation $\mu\omega_z =\omega_{\overline{z}}$.  

Set ${\cal QC}(G)$ to be the space of those representations
$\rho$ of $G$ into ${\rm PSL}_2({\bf C})$ which are induced by a
quasiconformal homeomorphism of $\overline{{\bf C}}$, so that $\rho\in
{\cal QC}(G)$ if there exists a quasiconformal homeomorphism
$\omega$ of $\overline{\bf C}$ so that $\rho(g)
=\omega\circ g\circ\omega^{-1}$ for all $g\in G$. By
definition, we have that ${\cal QC}(G)$ is contained in ${\cal
D}(G)$. Set ${\rm QC}(G)={\cal QC}(G)/{\rm PSL}_2({\bf C})$.

It is known that ${\rm QC}(G)$ is a complex manifold, and is actually
the quotient of the Teichm\"uller space of the (possibly disconnected)
quotient Riemann surface $\Omega(G)/G$ by a properly discontinuous
group of biholomorphic automorphims.  This result, in its full
generality, follows from the work of a number of authors, including
Maskit \cite{maskit-self-maps}, Kra \cite{kra-spaces}, Bers
\cite{bers-spaces}, and Sullivan \cite{sullivan-invariant}. 

We note here, in the case that $G$ is a geometrically finite,
minimally parabolic Kleinian group, that it follows from the
Isomorphism Theorem of Marden \cite{marden} that ${\rm QC}(G)$ is the
component of ${\rm MP}(G)$ containing the identity representation. 

Sullivan \cite{sullivan-invariant} has shown, for a finitely generated
Kleinian group $G$, if there exists a quasiconformal homeomorphism
$\omega$ of $\overline{\bf C}$ which conjugates $G$ to a Kleinian
group and which is conformal on $\Omega(G)$, then $\omega$ is
necessarily a M\"obius transformation.  In other words, if $\omega$
conjugates $G$ to subgroup of ${\rm PSL}_2({\bf C})$, then $\mu
=\omega_{\overline{z}}/\omega_z$ is equal to $0$ on $\Lambda(G)$.

In particular, if $\Omega(G)$ is empty, then ${\rm QC}(G)$ consists of
a single point, namely the identity representation.  This can be
viewed as a generalization of Theorem \ref{mostow-rigidity}, as
Sullivan's result also holds for an infinite volume hyperbolic
$3$--manifold $N$ whose uniformizing Kleinian group $G$ has limit set
the whole Riemann sphere.

We note here that the study of quasiconformal deformations of finitely
Kleinian groups is the origin of the Ahlfors Measure Conjecture.  In
\cite{ahlfors-finiteness}, Ahlfors raises the question of whether the
limit set of a finitely generated Kleinian group with non-empty domain
of discontinuity necessarily has zero area.  If this conjecture is
true, then it would be impossible for a quasiconformal deformation of
a finitely generated Kleinian group $G$ to be supported on the
limit set of $G$.  The result of Sullivan mentioned above implies
that no such deformation exists, though without solving the Measure 
Conjecture, which has not yet been completely resolved.  It is known
that the Measure Conjecture holds in a large number of cases, in
particular it holds for all topologically tame groups.  For a
discussion of this connection, we refer the reader to Canary
\cite{canary-tame} and the references contained therein.

There are several classes of Kleinian groups for which ${\rm
QC}(G)$ has been extensively studied, which we discuss here.

A {\em Schottky group} is a finitely generated, purely loxodromic
Kleinian group $G$ which is free on $g$ generators and whose
domain of  discontinuity is non-empty; the number of generators is
sometimes referred to as the {\em genus} of the Schottky group.  This
is not the original definition, but is equivalent to the usual
definition by a theorem of Maskit \cite{maskit-characterization}.  In
particular, a Schottky group is necessarily convex cocompact.
Chuckrow \cite{chuckrow} shows that any two Schottky groups of the
same rank are quasiconformally conjugate, so that ${\rm QC}(G)$
is in fact equal to the space ${\rm MP}(G)$ of all convex cocompact
representations of a group $G$ which is free on $g$ generators into
${\rm PSL}_2({\bf C})$.  

In the same paper \cite{chuckrow}, Chuckrow also engages in an
analysis of the closure of ${\rm QC}(G)$ in ${\rm R}(G)$ for
a Schottky group $G$ of genus $g$. In particular, she shows that
every point $\rho$ in $\partial {\rm QC}(G)$ has the property
that $\rho(G)$ is free on $g$ generators, and contains no
elliptic elements of infinite order.  However, this in itself is not
enough to show that $\rho(G)$ is discrete, as Greenberg
\cite{greenberg} has constructed free, purely loxodromic subgroups of
${\rm PSL}_2({\bf C})$ which are not discrete.  

More generally, Chuckrow also shows that the limit $\rho$ of a
convergent sequence $\{\rho_n\}$ of type-preserving faithful
representations in ${\cal HOM}(G)$ is again a faithful
representation of $G$, and that $\rho(G)$ contains no
elliptic elements of infinite order.

J\o rgensen \cite{jorgensen} credits his desire to generalize the
results of Chuckrow \cite{chuckrow} to leading him to what is now
commonly referred to as J\o rgensen's inequality, which states that if
$\gamma$ and $\varphi$ are elements of ${\rm PSL}_2({\bf C})$ which
generate a non-elementary Kleinian group, then $|{\rm tr}^2(\gamma)
-4| + |{\rm tr}([\gamma, \varphi] ) -2| \ge 1$, where ${\rm
tr}(\gamma)$ is the trace of a matrix representative of $\gamma$ in
${\rm SL}_2({\bf C})$.  The proof of Theorem
\ref{deformation-space-closed} is a direct application of this
inequality. 

For a Schottky group $G$, it is known that ${\rm AH}(G)$ is not
compact.  There is work of Canary \cite{canary-algebraic} and Otal
\cite{otal-degeneracy} on a conjecture of Thurston which gives
conditions under which sequences in ${\rm QC}(G)$ have convergent
subsequences; we do not give details here, instead referring the
interested reader to the papers cited above.

We also mention here the work of Keen and Series
\cite{keen-series-riley} on the Riley slice of the space of
$2$--generator Schottky groups, in which they introduce coordinates on
the Riley slice and study the cusp points on the boundary of the Riley
slice.

A {\em quasifuchsian} group is a finitely generated Kleinian group
whose limit set is a Jordan curve and which contains no element
interchanging the two components of its domain of discontinuity.
Consequently, every quasifuchsian group is isomorphic to the
fundamental group of a surface.  It is known that any two
isomorphic purely loxodromic quasifuchsian groups are quasiconformally
conjugate, by work of Maskit \cite{maskit-boundaries}, and hence
for a purely loxodromic quasifuchsian group $G$ we have that
${\rm MP}(G) ={\rm QC}(G)$.

This equality does not hold for quasifuchsian groups uniformizing
punctured surfaces, for several reasons.  First, the quasifuchsian
groups uniformizing the three-times punctured sphere and the
once-punctured torus are isomorphic, namely the free group of rank
two, but cannot be quasiconformally conjugate, as the surfaces are not
homeomorphic.  Second, as every quasifuchsian group isomorphic to the
free group $G$ of rank two contains parabolic elements, no
quasifuchsian group isomorphic to $G$ lies in ${\rm MP}(G)$.

It is known that ${\rm QC}(G)$ is biholomorphically equivalent to the
product of Teichm\"uller spaces ${\cal T}(S)\times {\cal
T}(\overline{S})$, where $S$ is one of the components of $\Omega(G)/G$
and $\overline{S}$ is its complex conjugate.

A {\em Bers slice} of ${\rm QC}(G)$ for a quasifuchsian group $G$ is a
subspace of ${\rm QC}(G)$ of the form $B(s_0) ={\cal T}(S)\times
\{s_0\}$.  The structure of the closure of $B(s_0)$ in ${\rm AH}(G)$
has been studied by a number of authors, including Bers
\cite{bers-boundaries}, Kerckhoff and Thurston
\cite{kerckhoff-thurston}, Maskit \cite{maskit-boundaries}, McMullen
\cite{mcmullen-cusps}, and Minsky \cite{minsky-punctured}.  In
particular, Bers \cite{bers-boundaries} showed that the closure
$\overline{B(s_0)}$ of $B(s_0)$ is compact, and Kerckhoff and Thurston
\cite{kerckhoff-thurston} have shown that the compactification
$\overline{B(s_0)}$ depends on the basepoint $s_0$, and so there are
actually uncountably many such compactifications.  Among
other major results, Minsky \cite{minsky-punctured} has shown that
every punctured torus group lies in the boundary of ${\rm QC}(G)$,
where $G$ is a quasifuchsian group uniformizing a punctured torus and
where a {\em punctured torus group} is a Kleinian group generated by
two elements with parabolic commutator.  In particular, this shows
that the relative version of the Density Conjecture holds for
punctured torus groups.

There are other slices of ${\rm QC}(G)$ which have been extensively
studied.  There is the extensive work of Keen and Series, see for
instance \cite{keen-series-bams}, \cite{keen-series-pleating}, and
\cite{keen-series-invariants}, inspired in part by unpublished work of
Wright \cite{wright}, on the Maskit slice of the Teichm\"uller 
space of a punctured torus in terms of {\em pleating coordinates},
which are natural and geometrically interesting coordinates on the
Teichm\"uller space of the punctured torus which are given in terms of
the geometry of the corresponding hyperbolic $3$--manifolds.

In the case that $G$ is a Kleinian group for which the corresponding
$3$--manifold $M =({\bf H}^3\cup\Omega(G))/G$ is a compact,
acylindrical $3$--manifold with non-empty, incompressible boundary,
then every representation in ${\rm MP}(G)$ in fact lies in ${\rm
QC}(G)$; this follows from work of Johannson \cite{johannson}.
In addition, Thurston \cite{thurston-acylindrical} has shown that
${\rm AH}(G)$ is compact for such $G$; another proof is given by
Morgan and Shalen \cite{ms-III}.  

\subsection{The spaces of ${\cal TT}(G)$ and ${\rm TT}(G) ={\cal
TT}(G)/{\rm PSL}_2({\bf C})$ of topologically tame representations} 

There is one last class of deformations which we need to define,
before beginning our discussion of the relationships between these
spaces.  We begin with a topological definition.  A compact
submanifold $M$ of a hyperbolic 3--manifold $N$ is a {\em compact core}
if the inclusion of $M$ into $N$ is a homotopy equivalence.  The Core
Theorem of Scott \cite{scott} implies that every hyperbolic
$3$--manifold with finitely generated fundamental group has a compact
core. Marden \cite{marden} asked whether every hyperbolic $3$--manifold
$N$ with finitely generated fundamental group is necessarily {\em
topologically tame}, in that $N$ is homeomorphic to the interior of
its compact core.

Set ${\cal TT}(G)$ to be the subspace of ${\cal D}(G)$ consisting of
the representations $\rho$ with minimally parabolic, topologically
tame image $\rho(G)$.\nl  Set ${\rm TT}(G) ={\cal TT}(G)/{\rm PSL}_2({\bf
C})$. 

There is a notion related to topological tameness, namely {\em
geometric tameness}, first defined by Thurston \cite{thurston-notes}.
We do not discuss geometric tameness here; the interested reader
should consult Thurston \cite{thurston-notes}, Bonahon \cite{bonahon},
or Canary \cite{canary-ends}.  Thurston \cite{thurston-notes} showed
that geometrically tame hyperbolic $3$--manifolds with freely
indecomposible fundamental group are topologically tame and satisfy
the Ahlfors Measure Conjecture.  Bonahon \cite{bonahon} showed that
if every non-trivial free product splitting of a finitely generated
Kleinian group $\Gamma$ has the property that there exists a parabolic
element of $\Gamma$ not conjugate into one of the free factors, then
$\Gamma$ is geometrically tame.  Canary \cite{canary-ends} extended
the definition of geometrically tame to all hyperbolic $3$--manifolds,
proved that topologically tame hyperbolic $3$--manifolds are
geometrically tame, and proved that topological tameness has a number
of geometric and analytic consequences; in particular, he established
that the Ahlfors Measure Conjecture holds for topologically tame
Kleinian groups. 

\section{Geometric limits}
\label{geometric-limits}

There is a second notion of convergence for Kleinian groups which is
distinct from the topology described above, which is equally important
in the study of deformations spaces.

A sequence $\{\Gamma_n\}$ of Kleinian groups converges {\em
geometrically} to a Kleinian group $\widehat\Gamma$ if two conditions
are met, namely that every element of $\widehat\Gamma$ is the limit of
a sequence of elements $\{\gamma_n\in\Gamma_n\}$ and that every
accumulation point of every sequence $\{\gamma_n\in\Gamma_n\}$ lies in
$\widehat\Gamma$.  Note that, unlike the topology of algebraic
convergence described above, the geometric limit of a sequence of
isomorphic Kleinian groups need not be isomorphic to the groups in the
sequence, and indeed need not be finitely generated.  However, it is
known that the geometric limit of a sequence of non-elementary,
torsion-free Kleinian groups is again torsion-free.

We note here that it is possible to phrase the definition of
geometric convergence in terms of the quotient hyperbolic
$3$--manifolds.  Setting notation, let $0$ denote a choice of basepoint
for ${\bf H}^3$, and let $p_j\co {\bf H}^3\to N_j ={\bf H}^3/\rho_j(G)$
and $p\co {\bf H}^3\rightarrow\widehat N ={\bf H}^3/\widehat\Gamma$ be
the covering maps. Let $B_{R}(0)\subset{\bf H}^3$ be a ball of radius
$R$ centered at the basepoint $0$.

\begin{lemma} A sequence of torsion-free Kleinian groups
$\{\Gamma_n\}$ converges geometrically to a torsion-free Kleinian
group $\widehat\Gamma$ if and only if there exists a sequence
$\{(R_n,K_n)\}$ and a sequence of orientation preserving maps
$\tilde{f_n} \co B_{R_n}(0) \to {\bf H}^3$ such that the following hold:
\begin{enumerate}
\item[\bf1\rm)] $R_n\to\infty$ and $K_n\to 1$ as $i\rightarrow\infty$;
\item[\bf2\rm)] the map $\tilde{f_n}$ is a $K_n$--bilipschitz diffeomorphism onto
its image, $\tilde{f_n}(0)=0$, and $\{\tilde{f_n} |_A\}$ converges to
the identity for any compact set $A$; and
\item[\bf3\rm)] $\tilde{f_n}$ descends to a map $f_n\co Z_n \to \widehat N$, where $Z_n
= B_{R_n}(0)/\Gamma_n$ is a submanifold of $N_n$; moreover, $f_n$ is
also an orientation preserving $K_n$--biLipschitz diffeomorphism onto
its image.
\end{enumerate}
\label{geom-convergence}
\end{lemma}

For a proof of this Lemma, see Theorem 3.2.9 of Canary, Epstein, and
Green \cite{canary-epstein-green}, and Theorem E.1.13 and Remark
E.1.19 of Benedetti and Petronio \cite{benedetti-petronio}.

A fundamental example of the difference between algebraic and
geometric convergence of Kleinian groups is given by the following
explicit example of J\o rg\-ensen and Marden \cite{jorgensen-marden};
earlier examples are given in J\o rgensen \cite{jorgensen-cyclic}.
Choose $\omega_1$ and $\omega_2$ in ${\bf C} -\{ 0\}$ which are
linearly independent over ${\bf R}$, and for each $n\ge 1$ set
$\omega_{1n} =\omega_1+n\omega_2$, $\omega_{2n} =\omega_2$, and
$\tau_n =\omega_{2n}/\omega_{1n}$. Consider the loxodromic elements
$L_n(z) =\exp(-2\pi {\bf i}\tau_n) z +\omega_2$.  Then, as
$n\rightarrow\infty$, $L_n$ converges to $L(z) =z+\omega_2$, and so
$\langle L_n\rangle$ converges algebraically to $\langle L\rangle$.
However, note that $L_n^{-n}(z)$ converges to $K(z)=z+\omega_1$ as
$n\rightarrow\infty$.  Hence, $\langle L_n\rangle$ converges
geometrically to $\langle L, K\rangle ={\bf Z}\oplus {\bf Z}$.

This example of the geometric convergence of loxodromic cyclic groups
to rank two parabolic groups underlies much of the algebra of the
operation of {\em Dehn surgery}, which we describe here.  

Let $M$ be a compact, hyperbolizable $3$--manifold, let $T$ be a
torus component of $\partial M$, and choose a meridian--longitude
system $(\alpha, \beta)$ on $T$.  Let $P$ be a solid torus and let $c$
be a simple closed curve on $\partial P$ bounding a disc in $P$.  For
each pair $(m,n)$ of relatively prime integers, let $M(m,n)$ be the
$3$--manifold by attaching $\partial P$ to $T$ by an
orientation-reversing homeomorphism which identifies $c$ with $m\alpha
+n\beta$; we refer to $M(m,n)$ as the result of $(m,n)$ Dehn surgery
along $T$.  The following Theorem describes the basic properties of
this operation; the version we state is due to Comar \cite{comar}.

\begin{theorem}[Comar \cite{comar}] Let $M$ be a compact,
hyperbolizable $3$--manifold and let $T =\{ T_1,\ldots, T_k\}$ be a
non-empty collection of tori in $\partial M$. Let $\hat N={\bf
H}^3/\Gamma$ be a geometrically finite hyperbolic $3$--manifold and let
$\psi\co {\rm int}(M)\rightarrow N$ be an orientation preserving
homeomorphism.  Further assume that every parabolic element of
$\Gamma$ lies in a rank two parabolic subgroup. Let $(m_i,l_i)$ be a
meridian--longitude basis for $T_i$. Let $\{ ({\bf p}_n, {\bf q}_n) =
((p^1_n,q^1_n), \ldots,(p^k_n,q^k_n)) \}$ be a sequence of $k$--tuples
of pairs of relatively prime integers such that, for each $i$, $\{
(p^i_n,q^i_n)\}$ converges to $\infty$ as $n\rightarrow\infty$. 

Then, for all sufficiently large $n$, there exists a representation
$\beta_n\co \Gamma\rightarrow {\rm PSL}_2({\bf C})$ with discrete image
such that 
\begin{enumerate}
\item[\bf1\rm)] $\beta_n(\Gamma)$ is geometrically finite, uniformizes
$M({\bf p}_n,{\bf q}_n)$, and every parabolic element of
$\beta_n(\Gamma)$ lies in a rank two parabolic subgroup; 
\item[\bf2\rm)] the kernel of $\beta_n\circ\psi_*$ is normally generated by $\{
m_1^{p_n^1} l_1^{q_n^1}, \ldots, m_k^{p_n^k} l_k^{q_n^k}\}$; and 
\item[\bf3\rm)] $\{ \beta_n\}$ converges to the identity representation of
$\Gamma$. 
\end{enumerate}
\label{dehn-surgery}
\end{theorem}

The idea of Theorem \ref{dehn-surgery} is due to Thurston
\cite{thurston-notes} in the case that the hyperbolic $3$--manifold $N$
has finite volume, so that $\partial M$ consists purely of tori.  In
this case, it is also known that ${\rm volume}({\bf
H}^3/\beta_n(\Gamma)) < {\rm volume}({\bf H}^3/\Gamma)$ for each $n$,
and that ${\rm volume}({\bf H}^3/\beta_n(\Gamma))\rightarrow {\rm
volume}({\bf H}^3/\Gamma)$ as $n\rightarrow\infty$.  For a more
detailed discussion of this phenomenon, we refer the reader to Gromov
\cite{gromov} and Benedetti and Petronio \cite{benedetti-petronio}.
The generalization to the case that $N$ has infinite volume is due
independently to Bonahon and Otal \cite{bonahon-otal} and Comar
\cite{comar}.  Note that the $\beta_n(\Gamma)$ are not isomorphic, and
hence there is no notion of algebraic convergence for these groups.

In the case that we have a sequence of representations in ${\cal
D}(G)$, the following result of J\o rgensen and Marden is extremely
useful. 

\begin{proposition}[J\o rgensen and Marden \cite{jorgensen-marden}]
Let $\{\rho_n\}$ be a sequence in ${\rm AH}(G)$ converging to $\rho$;
then, there is a subsequence of $\{\rho_n\}$, again called
$\{\rho_n\}$, so that $\{\rho_n(G)\}$ converges geometrically to
a Kleinian group $\widehat\Gamma$ containing $\rho(G)$. 
\label{alg-geom-conv}
\end{proposition}

A sequence $\{\rho_n\}$ in ${\cal D}(G)$ converges {\em strongly} to
$\rho$ if $\{\rho_n\}$ converges algebraically to $\rho$ and if
$\{\rho_n(G)\}$ converges geometrically to $\rho(G)$.  Note that we
may consider ${\cal D}(G)$, and ${\rm AH}(G)$, to be endowed with
topology of strong convergence, instead of the topology of algebraic
convergence.  We also refer the reader to the recent article of
McMullen \cite{mcmullen-dimI}, in which a variant of the notion of
strong convergence is explored in a somewhat more general setting.

Generalizing the behavior of the sequence of loxodromic cyclic groups
described above, examples of sequences $\{\rho_n\}$ in ${\cal D}(G)$
which converge algebraically to $\rho$ and for which $\{\rho_n(G)\}$
converges geometrically to a Kleinian group $\Gamma$ properly
containing $\rho(G)$ have been constructed by a number of authors,
including Thurston \cite{thurston-notes}, \cite{thurston-double},
Kerckhoff and Thurston \cite{kerckhoff-thurston}, Anderson and Canary
\cite{ac-books}, Ohshika \cite{ohshika-strong}, and Brock
\cite{brock-thesis}, \cite{brock-survey}, among others.  

J\o rgensen and Marden \cite{jorgensen-marden} carry out a very
detailed study of the relationship between the algebraic limit and the
geometric limit in the case when the geometric limit is assumed to be
geometrically finite.  In general, not much is known about the
relationship between the algebraic and geometric limits of a sequence
of isomorphic Kleinian groups.  We spend the remainder of this section
discussing this question.

A fundamental point in understanding how algebraic limits sit inside
geometric limits is the following algebraic fact, which is an easy
application of J\o rgensen's inequality.

\begin{proposition}[Anderson, Canary, Culler, and Shalen \cite{accs}]
Let $\{\rho_n\}$ be a sequence in ${\cal D}(G)$ which converges to
$\rho$ and for which $\{\rho_n(G)\}$ converges geometrically to a
Kleinian group $\widehat\Gamma$ containing $\rho(G)$.  Then, for each
$\gamma\in\widehat\Gamma -\rho(G)$, the intersection
$\gamma\rho(G)\gamma^{-1}$ is either trivial or parabolic cyclic.
\label{intersection}
\end{proposition}

One of the first applications of this result, also in \cite{accs}, was
to show, when the algebraic limit is a maximal cusp, that the convex
hull of the quotient $3$--manifold corresponding to the algebraic limit
embeds in the quotient $3$--manifold corresponding to the geometric
limit.  This was part of a more general attempt to understand the
relationship between the volume and the rank of homology for a finite
volume hyperbolic $3$--manifold.

Another application was given by Anderson and Canary
\cite{ac-coresI}. Before stating the generalization, we need to give a
definition.  Given a Kleinian group $\Gamma$, consider its associated
$3$--manifold $M =({\bf H}^3\cup\Omega(\Gamma))/\Gamma$, where
$\Omega(\Gamma)$ is the domain of discontinuity of $\Gamma$.  Then,
$\Gamma$ has connected limit set and no accidental parabolics if and
only if every closed curve $\gamma$ in $\partial M$ which is homotopic
to a curve of arbitrarily small length in the interior of $M$ with the
hyperbolic metric, is homotopic to a curve of arbitrarily small length
in $\partial M$, with its induced metric.

\begin{theorem}[Anderson and Canary \cite{ac-coresI}] Let $G$ be a
finitely generated, tor\-sion-free, non-abelian group, let $\{\rho_n\}$
be a sequence in ${\cal D}(G)$ converging to $\rho$, and suppose that
$\{\rho_n(G)\}$ converges geometrically to $\widehat\Gamma$.  Let $N
={\bf H}^3/\rho(G)$ and $\widehat N ={\bf H}^3/\widehat\Gamma$, and
let $\pi\co  N\rightarrow \widehat N$ be the covering map. If $\rho(G)$
has non-empty domain of discontinuity, connected limit set, and
contains no accidental parabolics, then there exists a compact core
$M$ for $N$ such that $\pi$ is an embedding restricted to $M$.
\label{ac-embedded-core}
\end{theorem}

One can apply Theorem \ref{ac-embedded-core} to show that certain
algebraically convergent sequences are actually strongly convergent.
This is of interest, as it is generally much more difficult to
determine strong convergence of a sequence of representations than to
determine algebraic convergence.

\begin{theorem}[Anderson and Canary \cite{ac-coresI}] Let $G$ be a
finitely generated, tor\-sion-free, non-abelian group and let
$\{\rho_n\}$ be a sequence in ${\cal D}(G)$ converging to $\rho$.
Suppose that $\rho_n(G)$ is purely loxodromic for all $n$, and that
$\rho(G)$ is purely loxodromic.  If $\Omega(\rho(G))$ is non-empty,
then $\{\rho_n(G)\}$ converges strongly to $\rho(G)$.  Moreover, 
$\{\Lambda(\rho_n(G))\}$ converges to $\Lambda(\rho(G))$.
\label{thm-ac-i}
\end{theorem}

\begin{theorem}[Anderson and Canary \cite{ac-coresI}] Let $G$ be a
finitely generated, tor\-sion-free, non-abelian group and let
$\{\rho_n\}$ be a sequence in ${\cal D}(G)$ converging to $\rho$.
Suppose that $\rho_n(G)$ is purely loxodromic for all $n$, that
$\rho(G)$ is purely loxodromic, and that $G$ is not a non-trivial free
product of (orientable) surface groups and cyclic groups, then
$\{\rho_n(G)\}$ converges strongly to $\rho(G)$.  Moreover,
$\{\Lambda(\rho_n(G))\}$ converges to $\Lambda(\rho(G))$. 
\label{thm-ac-ii}
\end{theorem}

Both Theorem \ref{thm-ac-i} and Theorem \ref{thm-ac-ii} have been
generalized by Anderson and Canary \cite{ac-coresII} to Kleinian
groups containing parabolic elements, under the hypothesis that the
sequences are type-preserving.

One reason that strong convergence is interesting is that strongly
convergent sequences of isomorphic Kleinian groups tend to be
extremely well behaved, as one has the geometric data coming from the
convergence of the quotient $3$--manifolds as well as the algebraic
data coming from the convergence of the representations.  For
instance, there is the following Theorem of Canary and Minsky
\cite{canary-minsky}.  We note that a similar result is proven
independently by Ohshika \cite{ohshika-limits}.  

\begin{theorem}[Canary and Minsky \cite{canary-minsky}] Let $M$ be a
compact, irreducible $3$--manifold and let $\{\rho_n\}$ be a sequence
in ${\rm TT}(\pi_1(M))$ converging strongly to $\rho$, where each
$\rho_n(\pi_1(M))$ and $\rho(\pi_1(M))$ are purely loxodromic.  Then,
$\rho(\pi_1(M))$ is 
topologically tame; moreover, for all sufficiently large $n$, there
exists a homeomorphism $\varphi_n\co  {\bf
H}^3/\rho_n(\pi_1(M))\rightarrow {\bf H}^3/\rho(\pi_1(M))$ so that
$(\varphi_n)_\ast =\rho\circ\rho_n^{-1}$.
\label{cm-strong}
\end{theorem}

By combining the results of Anderson and Canary \cite{ac-coresI} and of
Canary and Minsky \cite{canary-minsky} stated above, one may conclude
that certain algebraic limits of sequences of isomorphic topologically
tame Kleinian groups are again topologically tame.

There is also the following result of Taylor \cite{taylor}.

\begin{theorem}[Taylor \cite{taylor}] Let $G$ be a finitely
generated, torsion-free, non-abelian group, and let $\{\rho_n\}$ be a
sequence in ${\cal D}(G)$ converging strongly to $\rho$, where each
$\rho_n(G)$ has infinite co-volume.  If $\rho(G)$ is geometrically
finite, then $\rho_n(G)$ is geometrically finite for $n$ sufficiently
large.  
\end{theorem}

The guiding Conjecture in the study of the relationship between
algebraic and geometric limits, usually attributed to J\o rgensen, is
stated below.

\begin{conjecture}[J\o rgensen] Let $\Gamma$ be a finitely generated,
torsion-free, non-elementary Kleinian group, let $\{\rho_n\}$ be a
sequence in ${\cal D}(\Gamma)$ converging to $\rho$, and suppose that
$\{\rho_n(\Gamma)\}$ converges geometrically to $\widehat\Gamma$.  If
$\rho$ is type-preserving, then $\rho(\Gamma) =\widehat\Gamma$.  
\label{no-new-parabolics}
\end{conjecture}

As we have seen above, this conjecture has been shown to hold in a
wide variety of cases, including the case in which the sequence
$\{\rho_n\}$ is type-preserving and the limit group $\rho(\Gamma)$
either has non-empty domain of discontinuity or is not a non-trivial
free product of cyclic groups and the fundamental groups of orientable
surfaces.

\section{Functions on deformation spaces}
\label{functions-on-deformation}

There are several numerical quantities associated to a Kleinian group
$\Gamma$; one is the {\em Hausdorff dimension} ${\rm D}(\Gamma)$ of
the limit set $\Lambda(\Gamma)$ of $\Gamma$, another is the smallest
positive eigenvalue ${\rm L}(\Gamma)$ of the Laplacian on the
corresponding hyperbolic $3$--manifold ${\bf H}^3/\Gamma$.  These two
functions are closely related; namely, if $\Gamma$ is topologically
tame, then ${\rm L}(\Gamma) ={\rm D}(\Gamma)(2-{\rm D}(\Gamma))$ when
${\rm D}(\Gamma) \ge 1$, and ${\rm L}(\Gamma) =1$ when ${\rm
D}(\Gamma) \le 1$. The relationship between these two quantities has
been studied by a number of authors, including Sullivan
\cite{sullivan-entropy}, Bishop and Jones \cite{bishop-jones}, Canary
\cite{canary-laplacian}, and Canary, Minsky, and Taylor \cite{cmt}
(from which the statement given above is taken). It is natural to
consider how these functions behave on the spaces we have been
discussing in this note.

We begin by giving a few topological definitions.  A compact,
hyperbolizable $3$--manifold with incompressible boundary is a {\em
generalized book of $I$--bundles} if there exists a disjoint collection
$A$ of essential annuli in $M$ so that each component of the closure
of the complement of $A$ in $M$ is either a solid torus, a thickened
torus, or an $I$--bundle whose intersection with $\partial M$ is the
associated $\partial I$--bundle.

An {\em incompressible core} of a compact hyperbolizable $3$--manifold
is a compact submanifold $P$, possibly disconnected, with
incompressible boundary so that $M$ can be obtained from $P$ by adding
$1$--handles.

We begin with a pair of results of Canary, Minsky, and Taylor
\cite{cmt} which relates the topology of $M$ to the behavior of these
functions on a well-defined subset of ${\rm AH}(\pi_1(M))$, and show
that they are in a sense dual to one another.

\begin{theorem}[Canary, Minsky, and Taylor \cite{cmt}] Let $M$ be a
compact, hyperbolizable $3$--manifold.  Then, $\sup L(\rho(\pi_1(M)))
=1$ if and only if every component of the incompressible core of $M$
is a generalized book of $I$--bundles; otherwise, $\sup
L(\rho(\pi_1(M))) <1$. Here, the supremum is taken over all $\rho$ in
${\rm AH}(\pi_1(M))$ for which ${\bf H}^3/\rho(\pi_1(M))$ is
homeomorphic to the interior of $M$.
\end{theorem}

\begin{theorem}[Canary, Minsky, and Taylor \cite{cmt}] Let $M$ be a
compact, hyperbolizable $3$--manifold which is not a handlebody or a
thickened torus.  Then, $\inf D(\rho(\pi_1(M))) =1$ if and only if
every component of the incompressible core of $M$ is a geneneralized
book of $I$--bundles; otherwise, $\inf D(\rho(\pi_1(M))) >1$. Here, the
infimum is taken over all $\rho$ in ${\rm AH}(\pi_1(M))$ for which
${\bf H}^3/\rho(\pi_1(M))$ is homeomorphic to the interior of $M$. 
\end{theorem}

It is also possible to consider how these quantities behave under
taking limits.  We note that results similar to Theorem
\ref{hdim-continuous} have been obtained by McMullen
\cite{mcmullen-dimI}, who also shows that the function $D$ is not
continuous on ${\cal D}(\pi_1(M))$ in the case that $M$ is a
handlebody. 

\begin{theorem}[Canary and Taylor \cite{canary-taylor}] Let $M$ be a
compact, hyperbolizable $3$--manifold which is not homeomorphic to a
handlebody.  Then $D(\rho)$ is continuous on ${\cal D}(\pi_1(M))$
endowed with the topology of strong convergence.
\label{hdim-continuous}
\end{theorem}

Recently, Fan and Jorgenson \cite{fan-jorgenson} have made use of the
heat kernel to prove the continuity of small eigenvalues and small
eigenfunctions of the Laplacian for sequences of hyperbolic
$3$--manifolds converging to a geometrically finite limit manifold,
where the convergence is the variant of strong convergence considered 
by McMullen \cite{mcmullen-dimI}.

There are several functions on ${\rm QC}(G)$ which have been studied
by Bonahon.  In order to keep definitions to a minimum, we state his
results for geometrically finite $G$, though we note that they hold
for a general finitely generated Kleinian group $G$.  Given a
representation $\rho$ in ${\rm QC}(G)$, recall that the {\em convex
core} $C_\rho$ 
of ${\bf H}^3/\rho(G)$ is the smallest convex submanifold of ${\bf
H}^3/\rho(G)$ whose inclusion is a homotopy equivalence.  By
restricting the hyperbolic metric on ${\bf H}^3/\rho(G)$ to $\partial
C_\rho$, we obtain a map $\mu$ from ${\rm QC}(G)$ to the Teichm\"uller
space ${\cal T}(\Omega(G)/G)$ of the Riemann surface $\Omega(G)/G$.

\begin{theorem}[Bonahon \cite{bonahon-variations}] For a
geometrically finite Kleinian group $G$, the map $\mu\co  {\rm
QC}(G)\rightarrow {\cal T}(\Omega(G)/G)$ is continuously
differentiable. 
\end{theorem}

Another function on ${\rm QC}(G)$ studied by Bonahon, by developing an
analog of the Schl\"afli formula for the volume of a polyhedron in
hyperbolic space, is the function ${\rm vol}\co  {\rm QC}(G)\rightarrow
[0,\infty)$, which associates to $\rho\in {\rm QC}(G)$ the volume
${\rm vol}(\rho)$ of the convex core $C_\rho$ of ${\bf H}^3/\rho(G)$.

\begin{theorem}[Bonahon \cite{bonahon-schlafli}] Let $G$ be a
geometrically finite Kleinian group.  If the boundary $\partial
C_\rho$ of the convex core $C_\rho$ of ${\bf H}^3/\rho(G)$ is totally
geodesic, 
then $\rho$ is a local minimum of ${\rm vol}\co  {\rm QC}(G)\rightarrow
[0,\infty)$. 
\end{theorem}

It is known that the Hausdorff dimension of the limit set is a
continuous function on ${\rm QC}(\Gamma)$, using estimates relating
the Hausdorff dimension and quasiconformal dilitations due to Gehring
and V\"ais\"al\"a \cite{gehring-vaisala}.  In some cases, it is
possible to obtain more analytic information.

\begin{theorem}[Ruelle \cite{ruelle}] Let $\Gamma$ be a convex
cocompact Kleinian group whose limit set supports an expanding Markov
partition. Then, the Hausdorff dimension of the limit set is a real
analytic function on ${\rm QC}(\Gamma)$.
\label{expanding-markov}
\end{theorem}

Earlier work of Bowen \cite{bowen} shows that quasifuchsian and
Schottky groups support such Markov partitions.  The following
Theorem follows by combining these results of Bowen and Ruelle with a
condition which implies the existence of an expanding Markov
partition, namely that there exists a fundamental polyhedron in ${\bf
H}^3$ for the Kleinian group $G$ which has the {\em even cornered
property}, together with the Klein Combination Theorem.

\begin{theorem}[Anderson and Rocha \cite{anderson-rocha}] Let
$G$ be a convex cocompact\break Kleinian group which is isomorphic to the
free product of cyclic groups and fundamental groups of $2$--orbifolds.
Then, the Hausdorff dimension of the limit set is a real analytic
function on ${\rm QC}(G)$. 
\label{hdim-analytic}
\end{theorem}

We note here that it is not yet established that all convex cocompact
Kleinian groups support such Markov partitions.

Another function one can consider is the injectivity radius of the
corresponding quotient hyperbolic $3$--manifold.  For a hyperbolic
$3$--manifold $N$, the {\em injectivity radius} ${\rm inj}_N(x)$ at
a point $x\in N$ is one-half the length of the shortest homotopically
non-trivial closed curve through $x$.  The following Conjecture is due
to McMullen.

\begin{conjecture} Let $G$ be a finitely generated group with $g$
generators.  Then, there exists a constant $C =C(g)$ so that, if $N$
is a hyperbolic $3$--manifold with fundamental group isomorphic to $G$
and if $x$ lies in the convex core of $N$, then ${\rm inj}_N(x)\le C$.
\end{conjecture}

Kerckhoff and Thurston \cite{kerckhoff-thurston} show that, if $M$ is
the product of a closed, orientable surface $S$ of genus at least $2$
with the interval, then there exists a constant $C =C(M)$ so that if
$N$ is a hyperbolic $3$--manifold which is homeomorphic to the interior
of $M$ and if $N$ has no cusps, then the injectivity radius on the
convex core of $N$ is bounded above by $C$.  Fan \cite{fan}
generalizes this to show that, if $M$ is a compact, hyperbolizable 
$3$--manifold which is either a book of $I$--bundles or is acylindrical,
then there exists a constant $C =C(M)$ so that, if $N$ is any
hyperbolic $3$--manifold homeomorphic to the interior of $M$, then the
injectivity radius on the convex core of $N$ is bounded above by $C$.

We close by mentioning recent work of Basmajian and Wolpert
\cite{basmajian-wolpert} concerning the persistance of intersecting
closed geodesics.  Say that a Kleinian group $\Gamma$ has the {\em SPD
property} if all the closed geodesics in ${\bf H}^3/\Gamma$ are simple
and pairwise disjoint.

\begin{theorem}[Basmajian and Wolpert \cite{basmajian-wolpert}] Let
$G$ be a torsion-free, convex co-compact Kleinian group, and let
$U$ be the component of ${\rm MP}(G)$ containing the identity
representation.  Then, either
\begin{enumerate}
\item[\bf1\rm)] there exists a subset $V$ of $U$, which is the intersection of a
countably many open dense sets, so that $\rho(G)$ has the SPD property
for every $\rho\in V$, or 
\item[\bf2\rm)] there exists a pair of loxodromic elements $\alpha$ and $\beta$
of $G$ so that the closed geodesics in ${\bf H}^3/\rho(G)$
corresponding to loxodromic elements $\rho(\alpha)$ and $\rho(\beta)$
intersect at an angle constant over all $\rho\in U$; in particular,
there is {\rm no} element $\rho\in U$ so that $\rho(G)$ has the SPD
property. 
\end{enumerate}
\end{theorem}

They also show that the first possibility holds in the case that $G$
is a purely loxodromic Fuchsian group.

\bigskip
{\small \parskip 0pt \leftskip 0pt \rightskip 0pt plus 1fil \def\\{\par}
\sl\theaddress\par
\medskip
\rm Email:\stdspace\tt\theemail\par}
\recd

\end{document}